\documentclass[11pt]{article}

\usepackage{amsmath,amsthm,amsfonts,amssymb,latexsym}
\usepackage{amsmath}

\usepackage{a4wide}
\newtheorem{thm}{Theorem}[section]

\newtheorem{lem}[thm]{Lemma}
\newtheorem{prop}[thm]{Proposition}
\newtheorem{ex}[thm]{Example}
\newtheorem{rem}[thm]{Remark}

\newtheorem{theorem}{Theorem}[section]
\newtheorem{definition}[theorem]{Definition}

\newtheorem{proposition}[theorem]{Proposition}

\def\g{\mathfrak g}

\def\g{\mathfrak{g}}

  \def\Id{{\mathbb I}}

\def\L{{\cal   L}}

\def\L1#1{L^1(#1)}

\def\lef({\left(}
\def\rig){\right)}

\begin{document}


\date{}\title{Deformation of $\mathfrak{vect}(1)$-Modules of Symbols}

\author{ Imed Basdouri\thanks{D\'epartement de Math\'ematiques,
Facult\'e des Sciences de Sfax, BP 802, 3038 Sfax, Tunisie.
E.mail:basdourimed@yahoo.fr}\addtocounter{footnote}{0} \and
Mabrouk Ben Ammar
\thanks{D\'epartement de Math\'ematiques, Facult\'e des Sciences
de Sfax, BP 802, 3038 Sfax, Tunisie. E.mail:
mabrouk.benammar@fss.rnu.tn}\and B\'echir Dali
\thanks{D\'epartement de Math\'ematiques, Facult\'e des Sciences
de Bizerte, 7021 Zarzouna, Bizerte, Tunisie. E.mail:
bechir.dali@fss.rnu.tn} \and Salem Omri
\thanks{D\'epartement de Math\'ematiques, Facult\'e des Sciences
de Gafsa, Zarroug 2112, Tunisie. E.mail: omri\_salem@yahoo.fr} }

\maketitle

\begin{abstract}
We consider the action of the Lie algebra of polynomial vector
fields, $\mathfrak{vect}(1)$,  by the Lie derivative on
the space of symbols $\mathcal{S}_\delta^n=\bigoplus_{j=0}^n
\mathcal{F}_{\delta-j}$. We study deformations of this action. We
exhibit explicit expressions of some 2-cocycles generating  the
second cohomology space $\mathrm{H}^2_{\rm diff}(\mathfrak{vect}(1),{\cal D}_{\nu,\mu})$ where ${\cal
D}_{\nu,\mu}$ is the space of differential operators from
$\mathcal{F}_\nu$ to $\mathcal{F}_\mu$. Necessary second-order
integrability conditions of any infinitesimal deformations of
$\mathcal{S}_\delta^n$ are given. We describe completely the
formal deformations for some spaces $\mathcal{S}_\delta^n$ and we
give concrete examples of non trivial deformations.
\end{abstract}

\vfill\eject

\section{Introduction}
Let $\mathfrak{vect}(1)$ be the Lie algebra of polynomial
vector fields on $\mathbb{R}$. Consider the 1-parameter deformation
of the $\mathfrak{vect}(1)$-action on the space
$\mathbb{R}[x]$ of polynomial functions on $\mathbb{R}$ defined
by:\begin{equation*} L_{X\frac{d}{dx}}^\lambda(f)= Xf'+\lambda
X'f,\end{equation*} where $X, f\in\mathbb{R}[x]$ and
$X':=\frac{dX}{dx}$. Denote by $\cal{ F}_\lambda$ the
$\mathfrak{vect}(1)$-module structure on
$\mathbb{R}[x]$ defined by this action for a fixed $\lambda$.
Geometrically, $\cal{ F}_\lambda$ is the space of polynomial
weighted densities of weight $\lambda$ on $\mathbb{R}$:
\begin{equation*}{\cal F}_\lambda=\left\{ fdx^{\lambda}\mid f\in
\mathbb{R}[x]\right\}.
\end{equation*}
The space ${\cal F}_\lambda$ coincides with the space of vector
fields, functions and differential 1-forms for $\lambda = -1,\, 0$
and $1$, respectively.

Denote by ${\cal D}_{\nu,\mu}:=\mathrm{Hom}_{\text{diff}}({\cal
F}_\nu, {\cal F}_\mu)$ the $\mathfrak{vect}(1)$-module of
linear differential operators with the $\mathfrak{vect}(1)$-action given by the formula
\begin{equation}\label{Lieder2}L_X^{\nu,\mu}(A)=L_X^\mu\circ
A-A\circ L_X^\nu.\end{equation} Each module $\mathcal{D}_{\nu,\mu}$
has a natural filtration by the order of differential operators; the
graded module ${\cal S}_{\nu,\mu}:=\mathrm{gr}
\mathcal{D}_{\nu,\mu}$ is called the {\it space of symbols}. The
quotient-module
$\mathcal{D}^k_{\nu,\mu}/\mathcal{D}^{k-1}_{\nu,\mu}$ is isomorphic
to the module of weighted  densities $\mathcal{F}_{\mu-\nu-k} $, the
isomorphism is provided by the principal symbol map $\sigma_{\rm
pr}$ defined by:
\begin{equation*}
A=\sum_{i=0}^ka_i(x)\left(\frac{\partial}{\partial
x}\right)^i\mapsto\sigma_{\rm pr}(A)=a_k(x)(dx)^{\mu-\nu-k},
\end{equation*} (see, e.g., \cite{gmo}). As $\mathfrak{vect}(1)$-module, the space ${\cal S}_{\nu,\mu}$
depends only on the difference $\delta=\mu-\nu$, so that ${\cal
S}_{\nu,\mu}$ can be written as ${\cal S}_{\delta}$, and we have
$${\cal
S}_{\delta} = \bigoplus_{k=0}^\infty \mathcal{F}_{\delta-k}$$ as
$\mathfrak{vect}(1)$-modules. The space of symbols of
order $\leq n$ is $${\cal S}_\delta^n:=\bigoplus_{j=0}^n{\cal
F}_{\delta-j}.$$

The spaces ${\cal D}_{\nu,\mu}$ and $\mathcal{S}_{\delta}$ are not
isomorphic as $\mathfrak{vect}(1)$-modules:
$\mathcal{D}_{\nu,\mu}$ is a deformation of $\mathcal{S}_{\delta}$
in the sense of Richardson-Neijenhuis \cite{nr2}. In the last two
decades, deformations of various types of structures have assumed an
ever increasing role in mathematics and physics. For each such
deformation problem a goal is to determine if all related
deformation obstructions vanish and many beautiful techniques been
developed to determine when this is so. Deformations of Lie algebras
with base and versal deformations were already considered by
Fialowski in 1986 \cite{f1}. It was further developed, introducing a
complete local algebra base (local means a commutative algebra which
has a unique maximal ideal) by Fialowski in (1988) \cite{f2}. Also,
in \cite{f2}, the notion of miniversal (or formal versal)
deformation was introduced in general, and it was proved that under
some cohomology restrictions, a versal deformation exists. Later
Fialowski and Fuchs, using this framework, gave a construction for
versal deformation \cite{ff2}.

\vskip.3cm We use the framework of Fialowski \cite{f2} (see also
\cite{abbo} and \cite{aalo2}) and consider (multi-parameter)
deformations over complete local algebras. We construct the
miniversal deformation of this action and define the complete local
algebra related to this deformation.

According to Nijenhuis-Richardson \cite{nr2}, deformation theory of
modules is closely related to the computation of cohomology. More
precisely, given a Lie algebra $\frak{g}$ and a $\frak{g}$-module
$V$, the infinitesimal deformations of the $\frak{g}$-module
structure on $V$ , i.e., deformations that are linear in the
parameter of deformation, are related to
$\mathrm{H}^1\left(\mathfrak{g},\mathrm{End}(V)\right)$ .

Denote $\mathcal{D}:=\mathcal{D}(n,\delta)$ the $\mathfrak{vect}(1)$-module of differential operators on ${\cal
S}_\delta^n$. The infinitesimal  deformations of the $\mathfrak{vect}(1)$-module ${\mathcal S}^n_\beta$ are
classified  by the space
\begin{equation}\label{1}
\mathrm{H}^1_{\rm diff}\left(\mathfrak{vect}(1),
\mathcal{D}\right)=\oplus_{\lambda,k}\mathrm{H}^1_{\rm
diff}\left(\mathfrak{vect}(1),
\mathcal{D}_{\lambda,\lambda+k}\right),
\end{equation}
where  $\mathrm{H}^i_\mathrm{diff}$ denotes the differential
cohomology; that is, only cochains given by differential operators
are considered. Feigin and Fuchs computed  $\mathrm{H}^1_{\rm
diff}\left(\mathfrak{vect}(1),
\mathcal{D}_{\lambda,\lambda'}\right)$, see \cite{ff}. They showed
that non-zero cohomology $\mathrm{H}^1_{\rm diff}\left(\mathfrak{vect}(1),\mathcal{D}_{\lambda,\lambda'}\right)$ only
appear for particular values of weights that we call {\it resonant}
which satisfy  $\lambda'-\lambda\in\mathbb{N}$. Therefore, in
formula (\ref{1}), the summation $\oplus_{\lambda,k}$ is over all
$\lambda$ and $k$ satisfying
$0\leq\beta-\lambda-k\leq\beta-\lambda\leq n$.

In this paper we study the deformations of the structure of $\mathfrak{vect}(1)$-module on the space of symbols
$\mathcal{S}_\delta^n$. We give the second-order integrability
conditions which are sufficient in some cases. We will use the
framework of Fialowski \cite{f1,f2}  and Fialowski-Fuchs \cite{ff}
(see also \cite{abbo} and \cite{aalo2}) and consider
(multi-parameter) deformations over complete local algebra base. For
some examples, we will construct the miniversal deformation of this
action and define the local algebra related to this deformation. The
space $\mathrm{H}^1_{\rm diff}(\mathfrak{vect}(1),{\cal
D}_{\lambda,\lambda+k})$ was calculated in \cite{ff}, and for space
$\mathrm{H}^2_{\rm diff}(\mathfrak{vect}(1),{\cal
D}_{\lambda,\lambda+k})$ we can deduce the dimension  from
\cite{ff}, see also \cite{b}. We give explicit expressions of some
2-cocycles that span  $\mathrm{H}^2(\mathfrak{vect}(1),{\cal D}_{\lambda,\lambda+k})$.

\section{Cohomology Spaces}
Let $\frak{g}$ be a Lie algebra acting on a space $V$. The space
of $n$-cochains of $\frak{g}$ with values in $V$ is the
$\frak{g}$-module $$C^n(\frak{g}, V ) :=
\mathrm{Hom}(\wedge^n(\frak{g}),V).$$ The {\it coboundary
operator} $ \partial^n: C^n(\frak{g}, V )\rightarrow
C^{n+1}(\frak{g}, V )$ is a $\frak{g}$-map satisfying
$\partial^n\circ\partial^{n-1}=0$. The kernel of $\partial^n$,
denoted $Z^n(\mathfrak{g},V)$, is the space of $n$-{\it cocycles},
among them, the elements in the range of $\partial^{n-1}$ are
called $n$-{\it coboundaries}. We denote $B^n(\mathfrak{g},V)$ the
space of $n$-coboundaries.

By definition, the $n^{th}$ cohomolgy space is the quotient space
$$
\mathrm{H}^n
(\mathfrak{g},V)=Z^n(\mathfrak{g},V)/B^n(\mathfrak{g},V).
$$
We will only need the formula of $\partial^n$ (which will be
simply denoted $\partial$) in degrees 0, 1 and 2: for $v \in
C^0(\frak{g}, V ) = V$, $\partial v(X) := Xv$, for $ b\in
C^1(\frak{g}, V )$,
$$
\partial
b(X, Y ) := Xb(Y)-Y b(X) -b([X, Y ])
$$
and for $\Omega\in C^2(\frak{g}, V )$,
$$
\partial\Omega(X,Y,Z):=X\Omega(Y,Z)-\Omega([X,Y],Z)+\circlearrowleft
(X,Y,Z)
$$
where $  \circlearrowleft (X, Y,Z)$ denotes the summands obtained
from the two written ones by the cyclic permutation of the symbols
$X,\, Y,\,Z$.
\subsection{The First Cohomology Space}

The first cohomology space $\mathrm{H}^1_{\rm diff}(\mathfrak{vect}(1),{\cal D}_{\lambda,\lambda+k})$ was
calculated by Feigin and Fuks in \cite{ff}. The result is as
follows
\begin{thm}
\label{th1} The space $\mathrm{H}^1_{\rm diff}(\mathfrak{vect}(1),{\cal D}_{\lambda,\lambda+k})$ has the
following structure:
\begin{equation}
\label{CohSpace2} \mathrm{H}^1_{\rm
diff}(\mathfrak{vect}(1),{\mathcal{D}}_{\lambda,\lambda+k})\simeq\left\{
\begin{array}{ll}
\mathbb{R}&\hbox{ if }~~ k=0,2,3,4 \hbox{ for all
}\lambda,\\[2pt] \mathbb{R}^2& \hbox{ if
}~~\lambda=0\hbox{ ~and~
}k=1 ,\\[2pt] \mathbb{R}&\hbox{ if }~~
\lambda=0  \hbox{ or } \lambda=-4\hbox{ ~and~
}k=5, \\[2pt]
\mathbb{R}& \hbox{ if }~~ \lambda=-\frac{5\pm \sqrt{19}}{2}\hbox{
~and~
}k =6,\\[2pt] 0 &\hbox{ otherwise. }
\end{array}
\right.
\end{equation}
\end{thm}
These cohomology spaces are spanned by the cohomology classes of
the 1-cocycles, $C_{\lambda,\lambda+k}:\mathfrak{vect}(1)\rightarrow\mathcal{D}_{\lambda,\lambda+k
}$, that are collected in the following table. We write, for
$X\frac{d}{dx}\in\mathfrak{vect}(1)$ and
$f{dx}^{\lambda}\in{\cal F}_\lambda$,
\begin{align*}\begin{array}{llll}
C_{\lambda,\lambda+k
}(X\frac{d}{dx})(f{dx}^{\lambda})=C_{\lambda,\lambda+k
}(X)(f){dx}^{\lambda+k}. \end{array}\end{align*} \vskip.3cm Table 1.
Cocycles that span $\mathrm{H}^1_{\rm
diff}(\mathfrak{vect}(1),{\cal
D}_{\lambda,\mu})$\vskip.2cm
\begin{tabular}{|l|}\hline
  $C_{\lambda,\lambda}(X)(f)=X'f$ \\\hline
 $C_{0,1}(X)(f)=X''f$ \\\hline
   ${\widetilde C}_{0,1}(X)(f)=(X'f)'$ \\\hline
   $C_{\lambda,\lambda+2}(X)(f)=X^{(3)}f+2X''f'$ \\\hline
  $C_{\lambda,\lambda+3}(X)(f)=X^{(3)}f'+X''f''$ \\ \hline
  $C_{\lambda,\lambda+4}(X)(f)=-\lambda
X^{(5)}f+X^{(4)}f'-6X^{(3)}f''-4X''f^{(3)}$
  \\\hline
$C_{0,5}(X)(f)=2X^{(5)}f'-5X^{(4)}f''+10X^{(3)}f^{(3)}+5X''f^{(4)}$\\\hline

$C_{-4,1}(X)(f)=12X^{(6)}f+22X^{(5)}f'+5X^{(4)}f''-10X^{(3)}f^{(3)}-5X''f^{(4)}$\\\hline
  $C_{a_i,a_i+6}(X)(f)=\alpha_i X^{(7)}f-\beta_i X^{(6)}f'-\gamma_i
X^{(5)}f''-
  5X^{(4)}f^{(3)}+5X^{(3)}f^{(4)}+2X''f^{(5)}$  \\ \hline
\end{tabular}
\vskip.3cm where\begin{equation*}
\begin{array}{llllllll}a_1=-\frac{5+ \sqrt{19}}{2},
&\alpha_1=-\frac{22+ 5\sqrt{19}}{4}, &\beta_1=\frac{31+
7\sqrt{19}}{2}, &\gamma_1=\frac{25+
7\sqrt{19}}{2},\\[2pt]
a_2=-\frac{5- \sqrt{19}}{2}, &\alpha_2=-\frac{22- 5\sqrt{19}}{4},
& \beta_2=\frac{31- 7\sqrt{19}}{2},&\gamma_2=\frac{25-
7\sqrt{19}}{2}.\end{array}\end{equation*}

The maps $C_{\lambda,\lambda+j}(X)$ are naturally extended  to
${\cal S}_\delta^n=\bigoplus_{j=0}^n{\cal F}_{\delta-j}$.
\subsection{The Second Cohomology Space}
Let $\g$ a Lie algebra and $V$ a $\g$-module,  the {\it
cup-product} defined, for arbitrary linear maps $C_1 ,\, C_2 :
\frak g\rightarrow \mathrm{End}(V)$, is defined by:
\begin{equation}
\label{maurrer cartan1}
\renewcommand{\arraystretch}{1.4}
\begin{array}{l}
{}[\![C_1 , C_2]\!] : \frak g \otimes \frak g \rightarrow
\mathrm{End}(V)\\ {}[\![C_1 , C_2]\!] (x , y) = [C_1(x) , C_2(y)]
+ [C_2(x) , C_1(y)].\end{array}
\end{equation}
Therefore, it is easy to check that for any two $1$-cocycles $C_1$
and $C_2 \in Z^1 (\frak g , \mathrm{End}(V))$, the bilinear map
$[\![C_1 , C_2]\!]$ is a $2$-cocycle. Moreover, if one of the
cocycles $C_1$ or $C_2$ is a 1-coboundary, then $[\![C_1 ,
C_2]\!]$ is a $2$-coboundary. Therefore, we naturally deduce that
the operation (\ref{maurrer cartan1}) defines a bilinear map:
\begin{equation}
\label{cup-product} \mathrm{H}^1 (\frak g ,\mathrm{End}(
V))\otimes \mathrm{H}^1 (\frak g , \mathrm{End}( V))\rightarrow
\mathrm{H}^2 (\frak g , \mathrm{End}( V)).
\end{equation}
Thus, we can deduce the expressions of some 2-cocycles by computing
the cup-products of 1-cocycles. That is especially important if we
know the dimension of $\mathrm{H}^2 (\frak g , \mathrm{End}( V))$.

\vskip.3cm The second cohomology space $\mathrm{H}^2_{\rm
diff}(\mathfrak{vect}(1),{\cal D}_{\lambda,\lambda+k})$
of $\mathfrak{vect}(1)$ can be deduced from the work of
Feigin-Fuks \cite{ff} (see also \cite{b}). The result is as
follows
\begin{thm}
\label{th2} The space $\mathrm{H}^2_{\rm diff}(\mathfrak{vect}(1),{\cal D}_{\lambda,\lambda+k})$ has the
following structure:
\begin{align*}
 \mathrm{H}^2_{\rm diff}(\mathfrak{vect}(1),{\cal
D}_{\lambda,\lambda+k})\simeq \left\{
\begin{array}{ll}
\mathbb{R} &\text{ if }\left\{
\begin{array}{ll}~k=1,\lambda=0,\cr
~k= 2,3,4,7,8,9,10,11 \hbox{ for all }\lambda,\cr
~k=12,13,14\hbox{ but }\lambda\hbox{ is either
}\frac{1-k}{2}\hbox{ or
}\frac{1-k}{2}\pm\frac{\sqrt{12k-23}}{2},\cr
\end{array}
\right.\\[2pt]
\mathbb{R}^2 &\hbox{ if } ~~~~~~k=5,~\lambda=0,-4 ,\hbox{ or }k=6,~\lambda=a_1,a_2,\\[2pt]
 0  &\hbox{otherwise}.
\end{array}
\right.
\end{align*}
\end{thm}
\vskip.3cm In the sequel, we consider some 2-cocycles
$\Omega_{\lambda,\lambda+k }:\mathfrak{vect}(1)\times\mathfrak{vect}(1)\rightarrow\mathcal{D}_{\lambda,\lambda+k}$.
For $X\frac{d}{dx},\,Y\frac{d}{dx}\in\mathfrak{vect}(1)$
and $f{dx}^{\lambda}\in{\cal F}_\lambda$, we write
\begin{align*}\begin{array}{llll}
\Omega_{\lambda,\lambda+k
}(X\frac{d}{dx},Y\frac{d}{dx})(f{dx}^{\lambda})=\Omega_{\lambda,\lambda+k
}(X,Y)(f){dx}^{\lambda+k}. \end{array}\end{align*}

We need the following two lemmas:
\begin{lem}\label{lemc}
Let $b_{\lambda,\lambda+2}\in C^1
(\mathfrak{vect}(1),{\mathcal{D}}_{\lambda,\lambda+2})$
defined as follows: for
$X\frac{d}{dx}\in\mathfrak{vect}(1)$ and
$fdx^\lambda\in\mathcal{F}_\lambda$
\begin{equation}\label{cobo1}
b_{\lambda,\lambda+2}(X)(f)=\sum_{j=0}^3{}\alpha_jX^{(3-j)}f^{(j)}\end{equation}
where the coefficients $\alpha_j$ are constant. Then the map
$\partial
b_{\lambda,\lambda+2}:\mathfrak{vect}(1)\times
\mathfrak{vect}(1)
\rightarrow{\mathcal{D}}_{\lambda,\lambda+2}$ is given by
\begin{equation}\label{cobo2}\begin{array}{ll}
\partial
b_{\lambda,\lambda+2}(X,Y)(f)&=-\lambda\alpha_3XY^{(4)}f-\lambda\alpha_2X^{'}Y^{(3)}f-
\alpha_3(3\lambda+1)XY^{(3)}f^{'}-\alpha_2(2\lambda+1)X^{'}Y^{''}f^{'}\\&-
\alpha_3(3\lambda+3)XY^{''}f^{''}-\alpha_3XY^{'}f^{(3)}-(X\leftrightarrow
Y).\end{array}
\end{equation}
\end{lem}
\begin{proofname}. straightforward computation.\hfill$\Box$ \end{proofname}
\begin{lem}\label{lemc1}
Any $\partial b_{0,5}\in C^2
(\mathfrak{vect}(1),{\mathcal{D}}_{0,5})$ has
the general  following form: for
$X\frac{d}{dx},\,Y\frac{d}{dx}\in\mathfrak{vect}(1)$
and $f\in\mathcal{F}_0$,
\begin{equation}\begin{array}{llll}
\label{cobo3}
\partial
b_{0,5}(X,Y)(f)&=\alpha_0(9X^{''}Y^{(5)}+5X^{(3)}Y^{(4)})f-
\alpha_6XY^{(6)}f^{'}-\alpha_5X^{'}Y^{(5)}f^{'}\\&+(\alpha_2+5\alpha_1-\alpha_4)X^{''}Y^{(4)}f^{'}
-6\alpha_2XY^{(5)}f^{''}-5\alpha_5X^{'}Y^{(4)}f^{''}\\&+(2\alpha_2+3\alpha_3-4\alpha_4)X^{''}Y^{(3)}f^{''}
-15\alpha_6XY^{(4)}f^{(3)}-10\alpha_5X^{'}Y^{(3)}f^{(3)}
\\&-20\alpha_6XY^{(3)}f^{(4)}-10\alpha_5X^{'}Y^{''}f^{(4)}
-15\alpha_6XY^{''}f^{(5)}+ \alpha_5XY^{'}f^{(6)}-(X\leftrightarrow
Y).\end{array}
\end{equation} where the coefficients $\alpha_j$
are constant and the map $ b_{0,5}:\mathrm{Vect}_\mathrm{Pol}
\mathfrak{vect}(1)
\rightarrow{\mathcal{D}}_{0,5}$ is given by
\begin{equation}\label{cobo2}
b_{0,5}(X)(f)=\sum_{j=0}^6{}\alpha_jX^{(6-j)}f^{(j)}.\end{equation}
\end{lem}
\begin{proofname}. straightforward computation.\hfill$\Box$ \end{proofname}
\begin{prop}\label{ch1} The cohomology spaces $\mathrm{H}^2_{\rm
diff}(\mathfrak{vect}(1),{\mathcal{D}}_{\lambda,\lambda+k})$,
for $k=1,\,2,\,3,\,4$, are spanned by the cohomology classes of the
nontrivial 2-cocycles $\Omega_{\lambda,\lambda+k}$ defined by
\begin{align*}\Omega_{0,1}(X,Y)(f)~~~&=(X'Y''-X''Y')f,\\
\Omega_{\lambda,\lambda+2}(X,Y)(f)&=(X^{(3)}Y'-X'Y^{(3)})f+2(X''Y'-X'Y'')f',\\
\Omega_{\lambda,\lambda+3}(X,Y)(f)&=(X''Y^{(3)}-X^{(3)}Y'')f+(X^{(3)}Y'-X'Y^{(3)})f',\\
\Omega_{\lambda,\lambda+4}(X,Y)(f)&=-\lambda
X'Y^{(5)}f+X'Y^{(4)}f'-6X'Y^{(3)}f''-4X'Y''f^{(3)}-(X\leftrightarrow
Y)
\end{align*}
\end{prop}
\begin{proofname}.
The map $\Omega_{0,1}$ is the cup-product of the 1-cocyles $C_{0,1}$
and $C_{1,1}$. So, the map $\Omega_{0,1}$ is a 2-cocycle, therefore,
we will need only to prove that it is nontrivial since the space
$\mathrm{H}^2_{\rm
diff}(\mathfrak{vect}(1),{\mathcal{D}}_{0,1})$
is one dimensional. Let $b_{0,1}\in C^1
(\mathfrak{vect}(1),{\mathcal{D}}_{0,1})$
defined by
$$ b_{0,1}(X,f)=\sum_{j=0}^2{}\alpha_jX^{(2-j)}f^{(j)},\quad \text{where}\quad\alpha_j\in\mathbb{R},\,
X\frac{d}{dx}\in\mathfrak{vect}(1),\,f\in\mathcal{F}_0=\mathbb{R}[x].$$
Thus, $$\partial
b_{0,1}(X,Y)f=\alpha_2(X^{''}Y-XY^{''})f^{'}+\alpha_2(X^{'}Y-XY^{'})f^{''}$$
and therefore, it is clear that
$$\Omega_{0,1}(X,Y)f\neq\partial b_{0,1}(X,Y)f,\quad\text{for all}\quad \alpha_0,\,\alpha_1,\,\alpha_2\in\mathbb{R}.$$

The map $\Omega_{\lambda,\lambda+2}$ is the cup-product
$[\![C_{\lambda,\lambda+2},C_{\lambda,\lambda}]\!].$ By Lemma
\ref{lemc}, it is easy to check that $\Omega_{\lambda,\lambda+2}$
is a nontrivial 2-cocycle.

Besides, by direct computation, as before, we show that the
cup-products
$[\![C_{\lambda,\lambda+3},\Omega_{\lambda,\lambda}]\!]$ and
$[\![C_{\lambda+4,\lambda+4},\Omega_{\lambda,\lambda+4}]\!]$ are
nontrivial 2-cocycles. So, the spaces $\mathrm{H}^2_{\rm
diff}(\mathfrak{vect}(1),{\mathcal{D}}_{\lambda,\lambda+3})$
and  $\mathrm{H}^2_{\rm
diff}(\mathfrak{vect}(1),{\mathcal{D}}_{\lambda,\lambda+4})$
can be spanned respectively by the cohomology classes of the
nontrivial 2-cocycles $\Omega_{\lambda,\lambda+3}$ and
$\Omega_{\lambda,\lambda+4}$ defined by
\begin{equation*}
\Omega_{\lambda,\lambda+3}=[\![C_{\lambda,\lambda+3},
C_{\lambda,\lambda}]\!]\quad\text{and}\quad\Omega_{\lambda,\lambda+4}=[\![C_{\lambda+4,\lambda+4},
C_{\lambda,\lambda+4}]\!].
\end{equation*}
\hfill$\Box$
\end{proofname}

Now, we consider the cohomology spaces $\mathrm{H}^2_{\rm
diff}(\mathfrak{vect}(1),{\mathcal{D}}_{\lambda,\lambda+k})$
for $k=5,\,6.$ These spaces are generically trivial, but, for $k=5$
and $\lambda=-4, \,0$ or $k=6$ and $\lambda=a_1, \,a_2$ (where
$a_1=-\frac{5+ \sqrt{19}}{2}$ and $a_2=-\frac{5- \sqrt{19}}{2}$),
they are two dimensional. In the following proposition we exhibit a
basis for each of them.
\begin{prop}\label{ch2}
The cohomology spaces $\mathrm{H}^2_{\rm
diff}(\mathfrak{vect}(1),{\mathcal{D}}_{0,5})$,
$\mathrm{H}^2_{\rm
diff}(\mathfrak{vect}(1),{\mathcal{D}}_{-4,1})$
and  $\mathrm{H}^2_{\rm
diff}(\mathfrak{vect}(1),{\mathcal{D}}_{a_i,a_i+6})$,
(i=1, 2), are respectively spanned by the cohomology classes of the
nontrivial following 2-cocycles:
\begin{align*}\Omega_{0,5}(X,Y)(f)~~~~&=(X^{(5)}Y''+X^{(4)}Y^{(3)})f+
4X^{(4)}Y''f'+3X^{(3)}Y''f''-(X\leftrightarrow Y),\\
\widetilde{\Omega}_{0,5}(X,Y)(f)~~~~&=2X^{(5)}Y'f'+5X'Y^{(4)}f''+10X^{(3)}Y'f^{(3)}
+ 5X''Y'f^{(4)}-(X\leftrightarrow Y),\\
\Omega_{-4,1}(X,Y)(f)~~&=2X^{(4)}Y''f'+3X^{(3)}Y''f''-(X\leftrightarrow
Y),\\
\widetilde{\Omega}_{-4,1}(X,Y)(f)~~&=12X^{(6)}Y'f+22X^{(5)}Y'f'+5X^{(4)}Y'f''+
10X'Y^{(3)}f^{(3)}+5X'Y''f^{(4)}
\\&~~-(X\leftrightarrow Y),\\
\Omega_{a_i,a_i+6}(X,Y)(f)&=(X^{(5)}Y''+X^{(4)}Y^{(3)})f'+3X^{(4)}Y''f''
+2X^{(3)}Y''f^{(3)}-(X\leftrightarrow
Y)\\
\widetilde{\Omega}_{a_i,a_i+6}(X,Y)(f)&=\alpha_iX^{(7)}Y'f+\beta_iX'Y^{(6)}f'
+\gamma_iX'Y^{(5)}f''+ 5X'Y^{(4)}f^{(3)}+5X^{(3)}Y'f^{(4)}\\
&~~+2X''Y'f^{(5)}-(X\leftrightarrow Y),~~i=1,2.
\end{align*}
\end{prop}
\begin{proofname}. The  2-cocycles
$\Omega_{0,5}$ and $\widetilde{\Omega}_{0,5}$ are
 defined as follows:
 $$\Omega_{0,5}=[\![C_{2,5},\,{C}_{0,2}]\!],\quad\text{and}\quad
\widetilde{\Omega}_{0,5}=[\![C_{5,5},\,C_{0,5}]\!].$$ By Lemma
\ref{lemc1}, it is easy to show that these 2-cocyles are
nontrivial. Indeed, for instance, compering the term in $f$ in
both the expressions of $\Omega_{0,5}$ and of $\partial b_{0,5}$
given in (\ref{cobo3}), we see obviously that $\Omega_{0,5}$ can
not be a coboundary.

Similarly, we show that the 2-cocycles
$\Omega_{-4,1}=[\![C_{-1,1},\,{C}_{-4,-1}]\!]$,
$\widetilde{\Omega}_{-4,1}=C_{1,1},\,{C}_{-4,1}]\!]$,
$\Omega_{a_i,a_i+6}=[\![C_{a_i+3,a_i+6},\,{C}_{a_i,a_i+3}]\!]$ and
$\widetilde{\Omega}_{a_i,a_i+6}=[\![C_{a_i+6,a_i+6},\,{C}_{a_i,a_i+6}]\!]$
are  nontrivial. \hfill$\Box$
\end{proofname}

Now, we give basis for the spaces $\mathrm{H}^2_{\rm
diff}(\mathfrak{vect}(1),{\mathcal{D}}_{\lambda,\lambda+7})$
when $\lambda\notin\{0,-6\},  $ and for the spaces
$\mathrm{H}^2_{\rm
diff}(\mathfrak{vect}(1),{\mathcal{D}}_{\lambda,\lambda+8})$
when $2\lambda\neq -7\pm\sqrt{39}$.
\begin{prop}\label{ch3}
The spaces $\mathrm{H}^2_{\rm
diff}(\mathfrak{vect}(1),{\mathcal{D}}_{\lambda,\lambda+7})$
where $\lambda\neq0,-6$ and the spaces $\mathrm{H}^2_{\rm
diff}(\mathfrak{vect}(1),{\mathcal{D}}_{\lambda,\lambda+8})$
where $2\lambda\neq -7\pm\sqrt{39}$ are respectively spanned  by the
cohomology classes of the nontrivial following 2-cocycles:
\begin{align*}\Omega_{\lambda,\lambda+7}&=[\![C_{\lambda+3,\lambda+7},\,{C}_{\lambda,\lambda+3}]\!],\\
\Omega_{\lambda,\lambda+8}&=[\![C_{\lambda+4,\lambda+8},\,{C}_{\lambda,\lambda+4}]\!]
\end{align*}
\end{prop}
\begin{proofname}. By direct computation, we show that the cup-product
$[\![C_{\lambda+3,\lambda+7},\,{C}_{\lambda,\lambda+3}]\!]$  is a
nontrivial 2-cocycle if and only if  $\lambda\neq0,-6$. Similarly,
the cup-product
$[\![C_{\lambda+4,\lambda+8},\,{C}_{\lambda,\lambda+4}]\!]$  is a
nontrivial 2-cocycle if and only if  $2\lambda\neq
-7\pm\sqrt{39}$. \hfill$\Box$
\end{proofname}

For $k\geq9$, there are only  few cases where we can exhibit
2-cocycles by computation of cup-products of 1-cocycles. Theses
2-cocycles generating the corresponding cohomology spaces are
collected in the following proposition.
\begin{prop}\label{ch4}
The cohomology spaces $\mathrm{H}^2_{\rm
diff}(\mathfrak{vect}(1),{\mathcal{D}}_{a_i,a_i+9})$,
$\mathrm{H}^2_{\rm
diff}(\mathfrak{vect}(1),{\mathcal{D}}_{a_i-3,a_i+6})$,
$\mathrm{H}^2_{\rm
diff}(\mathfrak{vect}(1),{\mathcal{D}}_{-8,1})$,
$\mathrm{H}^2_{\rm
diff}(\mathfrak{vect}(1),{\mathcal{D}}_{0,9})$,
$\mathrm{H}^2_{\rm
diff}(\mathfrak{vect}(1),{\mathcal{D}}_{-4,5})$,
$\mathrm{H}^2_{\rm
diff}(\mathfrak{vect}(1),{\mathcal{D}}_{a_i,a_i+10})$,
$\mathrm{H}^2_{\rm
diff}(\mathfrak{vect}(1),{\mathcal{D}}_{a_i-4,a_i+6})$
are respectively spanned  by the cohomology classes of the
nontrivial following 2-cocycles:
$$
\begin{array}{llllllllllllllll}
\Omega_{a_i,a_i+9}&=&[\![C_{a_i+6,a_i+9},\,{C}_{a_i,\lambda+6}]\!],
&\Omega_{a_i-3,a_i+6}&=&[\![C_{a_i,a_i+6},\,{C}_{a_i-3,a_i}]\!],\\
\Omega_{-8,1}&=&[\![C_{-4,1},\,{C}_{-8,-4}]\!],
&\Omega_{0,9}&=&[\![C_{5,8},\,{C}_{0,5}]\!],\\
\Omega_{-4,5}&=&[\![C_{1,5},\,{C}_{-4,1}]\!],
&\Omega_{a_i,a_i+10}&=&[\![C_{a_i+6,a_i+10},\,{C}_{a_i,a_i+6}]\!],\\
\Omega_{a_i-4,a_i+6}&=&[\![C_{a_i,a_i+6},\,{C}_{a_i-4,a_i}]\!].
\end{array}
$$
\end{prop}
Here we omit the explicit expressions of these last 2-cocycles as
they are too long. But, as before, by direct computation, we show
that they are nontrivial.
\section{The General Framework}
In this section we define deformations of Lie algebra homomorphisms
and introduce the notion of miniversal deformations over complete
local algebras. Deformation theory of Lie algebra homomorphisms was
first considered with only one-parameter of deformation \cite{ff2,
nr2, r}. Recently, deformations of Lie algebras with
multi-parameters were intensively studied ( see,  e.g., \cite{abbo,
aalo2,  ro12, ro22}). Here we give an outline of this theory.
\subsection{Infinitesimal deformations}

Let $\rho_0:\frak g\to{\rm End}(V)$ be an action of a Lie algebra
$\frak g$ on a vector space $V$. When studying deformations of the
$\frak g$-action $\rho_0$, one usually starts with infinitesimal
deformations: $$ \rho=\rho_0+t\,C, $$ where $C:\frak g\to{\rm
End}(V)$ is a linear map and $t$ is a formal parameter. The
homomorphism condition $$ [\rho(x),\rho(y)]=\rho([x,y]), $$ where
$x,y\in\frak g$, is satisfied in order 1 in $t$ if and only if $C$
is a 1-cocycle. Moreover, two infinitesimal deformations $
\rho=\rho_0+t\,C_1, $ and $ \rho=\rho_0+t\,C_2, $ are equivalents
if and only if $C_1-C_2$ is a coboundary:
\begin{equation*}(C_1-C_2)(x)=[\rho_0(x),A]:=\partial
A(x),
\end{equation*}
where $A\in{\rm End}(V)$ and $\partial$ stands for differential of
cochains on $\frak g$ with values in $\mathrm{End}(V)$.  So, the
space $\mathrm{H}^1(\frak g,{\rm End}(V))$ determines and
classifies the  infinitesimal deformations up to equivalence.
(see,~e.g., \cite{Fuc2, nr2}). If $\mathrm{H}^1(\frak g,{\rm
End}(V))$ is multi-dimensional, it is natural to consider
multi-parameter deformations. More precisely, if
$\mathrm{dim}{\mathrm H}^1(\frak g,{\rm End}(V))=m$, then choose
1-cocycles $C_1,\ldots,C_m$ representing a basis of ${\mathrm
H}^1(\frak g,{\rm End}(V))$ and consider the infinitesimal
deformation
\begin{equation}
\label{InDefp} \rho=\rho_0+\sum_{i=1}^m{}t_i\,C_i,
\end{equation}
with independent parameters $t_1,\ldots,t_m$.

In our study, an infinitesimal deformation of the $\mathfrak{vect}(1)$-action on ${\cal S}^n_{\delta}$  is of
the form
\begin{equation}
\label{InDef}{\cal L}_X=L_X+{\cal L}^{(1)}_X,
\end{equation} where $L_X$ is the Lie derivative of ${\cal
S}^n_{\delta}$ along the vector field $X\frac{d}{dx}$ defined by
(\ref{Lieder2}), and
\begin{equation}
\label{InfinDef2} {\cal
L}_X^{(1)}=\sum_\lambda\sum_{j=0}^6{}t_{\lambda,\lambda+j}\,C_{\lambda,\lambda+j}(X)+{\widetilde
t}_{0,1}{\widetilde C}_{0,1}(X),
\end{equation}
and where $t_{\lambda,\lambda+j}$ and ${\widetilde t}_{0,1}$ are
independent parameters, $\delta-\lambda\in\mathbb{N}$,
$\delta-n\leq\lambda,\lambda+j\leq\delta$ and the 1-cocycles
$C_{\lambda,\lambda+j}$ and ${\widetilde C}_{0,1}$ are defined in
Table 1. We mention here that the term ${\widetilde
t}_{0,1}{\widetilde C}_{0,1}(X)$ don't appear in the expression of
${\cal L}_X^{(1)}$  if  $\delta-n\notin\mathbb{Z}_-$ or
$\delta\notin\mathbb{N}^*.$

\subsection{Integrability conditions}

Consider the problem of integrability of infinitesimal
deformations. Starting with the infinitesimal deformation
(\ref{InDefp}), we look for a formal series
\begin{equation}
\label{BigDef2} \rho= \rho_0+\sum_{i=1}^m{}t_i\,C_i+
\sum_{i,j}{}t_it_j\,\rho^{(2)}_{ij}+\cdots,
\end{equation}
where the highest-order terms
$\rho^{(2)}_{ij},\rho^{(3)}_{ijk},\ldots$ are linear maps from
$\frak g$ to ${\rm End(V)}$ such that \begin{equation} \label{map}
\rho:\frak g\to{\rm End(V)}\otimes\mathbb{C}[[t_1,\ldots,t_m]]
\end{equation} satisfies the homomorphism condition in any order in
$t_1,\ldots,t_m$.

However, quite often the above problem has no solution. Following
\cite{f1, f2} and \cite{aalo2}, we will impose extra algebraic
relations on the parameters $t_1,\ldots,t_m$. Let ${\cal R}$ be an
ideal in $\mathbb{C}[[t_1,\ldots,t_m]]$ generated by some set of
relations, the quotient
\begin{equation}
\label{TrivAlg2} {\cal A}=\mathbb{C}[[t_1,\ldots,t_m]]/{\cal R}
\end{equation}
is a local algebra with unity, and one can speak about  deformations
with base~${\cal A}$, see \cite{f1, f2} for details. The map
(\ref{map}) sends $\frak g$ to ${\rm End}(V)\otimes{\cal A}$.

{\ex \label{Example}  Consider the ideal $\cal R$ generated by all
the quadratic monomials $t_it_j$. In this case
\begin{equation}
\label{InfAlg2} {\cal A}=\mathbb{C}\oplus\mathbb{C}^m
\end{equation}
and any deformation is of the form (\ref{InDefp}). In this case
any infinitesimal deformation becomes a deformation with the base
$\cal A$ since $t_it_j=0$ in $\cal A$, for all $i,j=1,\ldots,m$. }

Given an infinitesimal deformation (\ref{InDefp}), one can always
consider it as a deformation with base (\ref{InfAlg2}). Our aim is
to find $\cal A$ which is big as possible, or, equivalently, we
look for relations on $t_1,\ldots,t_m$ which are necessary and
sufficient for integrability ( cf.\cite{abbo}, \cite{aalo2}).

\subsection{Equivalence and the miniversal deformation}

The notion of equivalence of deformations over commutative
associative algebras has been considered in \cite{ff2}.

{\definition  Two deformations, $\rho$ and $\rho'$ with the same
base $\cal A$ are called equivalent if there exists an inner
automorphism $\Psi$ of the associative algebra ${\rm
End}(V)\otimes\cal A$ such that
$$
\Psi\circ\rho=\rho'\hbox{ and }\Psi(\mathbb{I})=\mathbb{I},
$$
where $\mathbb{I}$ is the unity of the algebra ${\rm
End(V)}\otimes\cal A$. }

The following notion of miniversal deformation is fundamental. It
assigns to a~$\frak g$-module $V$ a canonical commutative
associative algebra $\cal A$ and a canonical deformation with base
$\cal A$.

{\definition A deformation $\rho$ with base $\cal A$ is called
miniversal, if
\begin{itemize}
  \item [(i)] for any other deformation, $\rho'$ with base (local)
$\cal A'$, there exists a  homomorphism $\psi:{\cal A}'\to{\cal A}$
satisfying $\psi(1)=1$, such that
$$\rho=(\Id\otimes\psi)\circ\rho'.$$
  \item [(ii)] in the notations of (i), if $\mathcal{A}$ is infinitesimal
  then  $\psi$ is unique.
\end{itemize}
If $\rho$ satisfies only the condition (i), then it is called
versal. }

The miniversal deformation corresponds to the smallest ideal $\cal
R$. We refer to \cite{ff2} for a construction of miniversal
deformations of Lie algebras and to \cite{aalo2} for miniversal
deformations of $\frak g$-modules.

\section{Second-order Integrability Conditions}

Assume that the infinitesimal deformation (\ref{InDef}) can be
integrated to a formal deformation

\begin{equation}
\label{formal1} {\cal L}_X=L_X+{\cal L}^{(1)}_X+{\cal
L}^{(2)}_X+{\cal L}^{(3)}_X+\cdots
\end{equation}
where ${\cal L}^{(1)}_X$ is given by (\ref{InfinDef2}) and ${\cal
L}^{(2)}_X$ is a quadratic polynomial in the parameters
$t_{\lambda,\lambda+k}$ with coefficients in ${\cal
D}^n_{\delta}$. We compute the conditions for the second-order
terms~${\cal L}^{(2) }$. Consider the quadratic terms of the
homomorphism condition
\begin{equation}
\label{HomomCond2} [{\cal L}_X,{\cal L}_Y]={\cal L}_{[X,Y]}.
\end{equation}
The homomorphism condition (\ref{HomomCond2}) gives for the
second-order terms the following (Maurer-Cartan) equation

\begin{equation}\label{MC1}\partial{\cal L}^{(2)}=-\frac{1}{2}\,[\![{\cal
L}^{(1)},{\cal L}^{(1)}]\!], \end{equation} so that the right hand
side of  (\ref{MC1}) is automatically a 2-cocycle. In our case, we
obtain:\begin{equation}\label{cap}\partial{\cal
L}^{(2)}=-\frac{1}{2}[\![\sum_\lambda\sum_{j=0}^6{}t_{\lambda,
\lambda+j}\,C_{\lambda,\lambda+j}+{\widetilde t}_{0,1}{\widetilde
C}_{0,1},\sum_\lambda\sum_{j=0}^6{}t_{\lambda,\lambda+j}\,C_{\lambda,\lambda+j}+{\widetilde
t}_{0,1}{\widetilde C}_{0,1}]\!].\end{equation}

\bigskip
Let us consider the 2-cocycles: $B_{\lambda,\lambda+k}\in Z^2_{\rm
diff}(\mathfrak{vect}(1),{\cal D}_{\lambda,\lambda+k})$,
for $k=0,\dots,10,$ defined by:
$$
B_{\lambda,\lambda+k}=\sum_{j=2}^k{}t_{\lambda+j,\lambda+k}{}t_{\lambda,\lambda+j}
\,[\![C_{\lambda+j,\lambda+k},C_{\lambda,\lambda+j}]\!].
$$
(We consider also $\widetilde{C}_{0,1}$ in the expression of
$B_{\lambda,\lambda+k}(X,Y)$ when it is possible: $\lambda=0$ and
$j=1$ or $\lambda+j=0$ and $\lambda+k=1$). Necessary conditions
for the integrability of the infinitesimal deformation
(\ref{InfinDef2}) are that any 2-cocycle $B_{\lambda,\lambda+k}$,
$k=0,\dots,10,$ must be a coboundary:
\begin{equation}\label{cob}B_{\lambda,\lambda+k}=\partial
b_{\lambda,\lambda+k},\end{equation} where $
b_{\lambda,\lambda+k}\in C^1(\mathfrak{vect}(1),{\cal
D}_{\lambda,\lambda+k})$. We easily see that
$B_{\lambda,\lambda}=0$, so, there are no integrability conditions
for $k=0$. In the following, we study, successively, the
second-order integrability  conditions for $k=1,\dots,10$.

\bigskip
\begin{prop}\label{pr1}For $k=1,\,2,\,3$, we have the following
second-order integrability  conditions of the infinitesimal
deformation~(\ref{InDef}):
\begin{equation}\label{k123}\begin{array}{llll}
t_{0,1}(t_{0,0}-t_{1,1})-t_{1,1}{\widetilde t}_{0,1}&=&0,~~(k=1,\lambda=0)\\[2pt]
t_{\lambda,\lambda+2}(t_{\lambda,\lambda}-t_{\lambda+2,\lambda+2})&=&0,~~(k=2)\\[2pt]
t_{\lambda,\lambda+3}(t_{\lambda,\lambda}-t_{\lambda+3,\lambda+3})&=&0,~~(k=3).\\[2pt]
\end{array}\end{equation}

\end{prop}
\begin{proofname}. Obviously, for $k=1$ and for
$\lambda\neq0$, we have $B_{\lambda,\lambda+1}=0$, since,
$C_{\lambda,\lambda+1}=0$, therefore, there are no conditions in
this case. For $\lambda=0$, we have
$$B_{0,1}=t_{0,1}t_{0,0}[\![C_{0,1},C_{0,0}]\!]+t_{1,1}t_{0,1}[\![C_{1,1},C_{0,1}]\!]
+t_{1,1}{\widetilde t}_{0,1}[\![\widetilde{C}_{0,1},C_{0,0}]\!].$$
By a straightforward computation, we show that
$$
[\![C_{1,1},C_{0,1}]\!]=-[\![C_{0,1},C_{0,0}]\!]=[\![\widetilde{C}_{0,1},C_{0,0}]\!]=\Omega_{0,1}.
$$
Therefore
$$
B_{0,1}=
\begin{array}{ll}
(t_{1,1}t_{0,1}+t_{1,1}{\widetilde
t}_{0,1}-t_{0,1}t_{0,0}){}\Omega_{0,1}.
\end{array}
$$
According to Proposition \ref{ch1}, $\Omega_{0,1}$ is a nontrvial
2-cocycle, so, the first integrability condition:
$t_{1,1}t_{0,1}+t_{1,1}{\widetilde t}_{0,1}-t_{0,1}t_{0,0}=0$,
holds

For $k=2$, we have $$B_{\lambda,\lambda+2}=t_{\lambda,\lambda+2}
t_{\lambda,\lambda}[\![C_{\lambda,\lambda+2},C_{\lambda,\lambda}]\!]+t_{\lambda,\lambda+2}t_{\lambda+2,\lambda+2}
[\![C_{\lambda+2,\lambda+2},C_{\lambda,\lambda+2}]\!].
$$
But, it is easy to show that
$[\![C_{\lambda,\lambda+2},C_{\lambda,\lambda}]\!]=
-[\![C_{\lambda+2,\lambda+2},C_{\lambda,\lambda+2}]\!]=\Omega_{\lambda,\lambda+2}$.
Thus, we get the integrability condition: $t_{\lambda,\lambda+2}
(t_{\lambda,\lambda}-t_{\lambda+2,\lambda+2})=0$, since, by
Proposition \ref{ch1}, $\Omega_{\lambda,\lambda+2}$ is a
nontrivial 2-cocycle.

Similarly, for $k=3$, we obtain
$t_{\lambda,\lambda+3}(t_{\lambda,\lambda}-t_{\lambda+3,\lambda+3})=0$.
\hfill$\Box$
\end{proofname}

\begin{prop}\label{pr3}For $k=4,\,5,\,6$, we have the following
second-order integrability  conditions of the infinitesimal
deformation~(\ref{InDef}), where in the first line
$\lambda\notin\{0,-3\}$:
\begin{equation}\label{k456}\begin{array}{llll}
t_{\lambda,\lambda+4} (t_{\lambda,\lambda}-t_{\lambda+4,\lambda+4})&=&0,\\[2pt]
t_{-3,1}
(t_{-3,-3}-t_{1,1})-\frac{1}{10}t_{-3,0}\widetilde{t}_{0,1}&=&0,\\[2pt]
t_{0,4}
(t_{0,0}-t_{4,4})+\frac{1}{10}\widetilde{t}_{0,1}t_{1,4}&=&0,\\[2pt]
30t_{0,0}t_{0,5}-12t_{0,1}t_{1,5}-12\widetilde{t}_{0,1}t_{1,5}+t_{0,2}t_{2,5}&=&0,\\[2pt]
t_{0,0}t_{0,5}-\frac{2}{5}\widetilde{t}_{0,1}t_{1,5}-t_{0,5}t_{5,5}&=&0,\\[2pt]
30t_{-4,-4}t_{-4,1}-12t_{-4,0}t_{0,1}-12t_{-4,0}\widetilde{t}_{0,1}
+t_{-4,-1}t_{-1,1}&=&0,\\[2pt]
t_{-4,1}t_{1,1}+\frac{2}{5}t_{-4,0}
\widetilde{t}_{0,1}-t_{-4,-4}t_{-4,1}&=&0,\\[2pt]
t_{a_i,a_i}t_{a_i,a_i+6}-t_{a_i,a_i+6}t_{a_i+6,a_i+6}&=&0,\\[2pt]
t_{a_i+3,a_i+6}t_{a_i,a_i+3}-R_it_{a_i,a_i}t_{a_i,a_i+6} +S_i
t_{a_i,a_i+2}t_{a_i+2,a_i+6}-T_it_{a_i,a_i+4}t_{a_i+4,a_i+6}&=&0.
\end{array}\end{equation}

where $R_i$, $S_i$ and $T_i$ are some constants which we leave out
their explicit expressions as they are so complicated.
\end{prop}
\begin{proofname}. 1) First, we show that
$[\![C_{\lambda+2,\lambda+4},C_{\lambda,\lambda+2}]\!]=\partial
b_{\lambda,\lambda+k}$ where
\begin{align*}\begin{array}{l} \label{}b_{\lambda,\lambda+4}(X)f=
\frac{2}{5}(\lambda-1)X^{(5)}f-2X^{(4)}f'.
\end{array}\end{align*}
For $\lambda\notin\{0,-3\}$, the 2-cocycle $B_{\lambda,\lambda+4}$
is defined by
\begin{align*}
B_{\lambda,\lambda+4}=\,\,&t_{\lambda,\lambda+4}
t_{\lambda,\lambda}[\![C_{\lambda,\lambda+4},C_{\lambda,\lambda}]\!]+t_{\lambda,\lambda+4}
t_{\lambda+4,\lambda+4}[\![C_{\lambda+4,\lambda+4},C_{\lambda,\lambda+4}]\!]\\&+t_{\lambda,\lambda+2}
t_{\lambda+2,\lambda+4}[\![C_{\lambda+2,\lambda+4},C_{\lambda,\lambda+2}]\!].
\end{align*}
We check that $[\![C_{\lambda,\lambda+4},C_{\lambda,\lambda}]\!]=
-[\![C_{\lambda+4,\lambda+4},C_{\lambda,\lambda+4}]\!]=\Omega_{\lambda,\lambda+4}$.
Thus, we get the following condition: $$t_{\lambda,\lambda+4}
(t_{\lambda,\lambda}-t_{\lambda+4,\lambda+4})=0.$$ For
$\lambda=0$, the cup-products $[\![C_{1,4},C_{0,1}]\!]$ and
$[\![C_{1,4},\widetilde{C}_{0,1}]\!]$ appear in the expression of
$B_{0,4}$. But, we check that
$[\![C_{1,4},\widetilde{C}_{0,1}]\!]={1\over10}\Omega_{0,4}$ and
$[\![C_{1,4},\widetilde{C}_{0,1}]\!]=\partial \widetilde{b}_{0,4}$
where \begin{align*}\begin{array}{l}
\label{}\widetilde{b}_{0,4}(X)f=- \frac{1}{5}X^{(5)}f+X''f^{(3)}.
\end{array}\end{align*}
For $\lambda=-3$, the cup-products $[\![C_{0,1},C_{-3,0}]\!]$ and
$[\![\widetilde{C}_{0,1},C_{-3,0}]\!]$ appear in the expression of
$B_{-3,1}$. We check that
$[\![\widetilde{C}_{0,1},C_{-3,0}]\!]=-{1\over10}\Omega_{-3,1}$ and
$[\![C_{1,4},\widetilde{C}_{0,1}]\!]=\partial \widetilde{b}_{-3,1}$
where \begin{align*}\begin{array}{l}
\label{}\widetilde{b}_{-3,1}(X)f=-
\frac{3}{10}X^{(5)}f-{1\over2}X^{(4)}f'.
\end{array}\end{align*}
Thus, for $\lambda\in\{0,\,-3\}$, we obtain the following
integrability conditions:
\begin{align*}\begin{array}{l}
t_{-3,1}
(t_{-3,-3}-t_{1,1})-\frac{1}{10}t_{-3,0}\widetilde{t}_{0,1}=t_{0,4}
(t_{0,0}-t_{4,4})+\frac{1}{10}\widetilde{t}_{0,1}t_{1,4}=0.
\end{array}\end{align*}

2) Now, if $\lambda\notin\{0,-4\}$  the 2-cocycle
$B_{\lambda,\lambda+5}$ is defined by
\begin{align*}
B_{\lambda,\lambda+5}=t_{\lambda,\lambda+2}
t_{\lambda+2,\lambda+5}[\![C_{\lambda+2,\lambda+5},C_{\lambda,\lambda+2}]\!]+t_{\lambda,\lambda+3}
t_{\lambda+3,\lambda+5}[\![C_{\lambda+3,\lambda+5},C_{\lambda,\lambda+3}]\!].
\end{align*} But, by a direct computation, we show that
$$
[\![C_{\lambda+2,\lambda+5},C_{\lambda,\lambda+2}]\!]=\partial
b_{\lambda,\lambda+5}\quad\text{and}\quad
[\![C_{\lambda+3,\lambda+5},C_{\lambda,\lambda+3}]\!]=\partial
\widetilde{b}_{\lambda,\lambda+5}
$$
where, for $\lambda\neq-2$,
\begin{align*}\begin{array}{llll}b_{\lambda,\lambda+5}(X)(f)&=\frac{-1}{\lambda^2+6\lambda+8}
\Big(\frac{7}{30}(\lambda+4)
X^{(6)}f+\frac{1}{10\lambda}(10\lambda^2+39\lambda-4)
X^{(5)}f'\\[2pt]
&~+\frac{1}{2\lambda}(4\lambda^2+17\lambda+4) X^{(4)}f'' +
\frac{1}{3\lambda}(3\lambda^2+11\lambda-4)X^{(3)}f^{(3)}\Big),
\\[2pt]\widetilde{b}_{\lambda,\lambda+5}(X)(f)&=\frac{-1}{\lambda^2+6\lambda+8}
\Big(-\frac{7}{30} \lambda X^{(6)}f
+\frac{21}{10}X^{(5)}f'+(\lambda+\frac{11}{2}) X^{(4)}f'' +
(\lambda+\frac{13}{3}) X^{(3)}f^{(3)}\Big)
\end{array}\end{align*}and
\begin{align*}\begin{array}{llll}b_{-2,3}(X)(f)&=-\frac{1}{12} X^{(5)}f'
-\frac{1}{3}X^{(4)}f'' -\frac{2}{3}X^{(3)}f^{(3)}
-\frac{7}{12}X''f^{(4)}\\[2pt]
\widetilde{b}_{-2,3}(X)(f)&=-\frac{7}{12} X^{(5)}f'
-\frac{5}{3}X^{(4)}f''-\frac{5}{3}X^{(3)}f^{(3)} -X''f^{(4)}
\end{array}\end{align*}

For $\lambda=0$, recall that the cohomology space $\mathrm{H}^2_{\rm
diff}(\mathfrak{vect}(1),{\mathcal{D}}_{0,5})$
is spanned by the cohomology classes of the 2-cocycles
$\Omega_{0,5}=[\![C_{2,5},\,{C}_{0,2}]\!]$ and $
\widetilde{\Omega}_{0,5}=[\![C_{5,5},\,C_{0,5}]\!]$ (see Proposition
\ref{ch2}). Moreover, in this case, we have
\begin{align*}
B_{0,5}=t_{0,1} t_{1,5}[\![C_{1,5},C_{0,1}]\!]+\widetilde{t}_{0,1}
t_{1,5}[\![C_{1,5},\widetilde{C}_{0,1}]\!]+t_{0,0}
t_{0,5}[\![C_{0,5},C_{0,0}]\!].
\end{align*}
But, by direct computation, we check that
\begin{align*}\begin{array}{llllllllll}
[\![C_{1,5},\,{C}_{0,1}]\!]&=-12\Omega_{0,5}+\partial
b_{0,5},\\[2pt]
[\![C_{1,5},\,\widetilde{C}_{0,1}]\!]&=-12\Omega_{0,5}+\frac{2}{5}\widetilde{\Omega}_{0,5}+\partial
\widetilde{b}_{0,5},\\[2pt]
[\![C_{0,5},\,{C}_{0,0}]\!]&=30\,\Omega_{0,5}-\widetilde{\Omega}_{0,5}+\partial
\overline{b}_{0,5}
\end{array}\end{align*} where
\begin{align*}\begin{array}{llll}b_{0,5}(X)(f)&=-
X^{(6)}f-55X^{(4)}f''-20X'f^{(5)},\\
\widetilde{b}_{0,5}(X)(f)&=- X^{(6)}f-45X^{(4)}f''-15X''f^{(4)}+
\frac{1}{5}X'f^{(5)},\\\overline{b}_{0,5}(X)(f)&=3 X^{(6)}f
+135X^{(4)}f''+45X''f^{(4)}.
\end{array}\end{align*}

Similarly, for $\lambda=-4$, the cohomology space $\mathrm{H}^2_{\rm
diff}(\mathfrak{vect}(1),{\mathcal{D}}_{-4,1})$
is spanned by the cohomology classes of
$\Omega_{-4,1}=[\![C_{-1,1},\,{C}_{-4,-1}]\!]$ and
$\widetilde{\Omega}_{-4,1}=C_{1,1},\,{C}_{-4,1}]\!]$, and we check
that
\begin{align*}\begin{array}{llllllllll}
[\![\widetilde{C}_{0,1},\,{C}_{-4,0}]\!]&=-12\Omega_{-4,1}+\frac{2}{5}\widetilde{\Omega}_{-4,1}+\partial
\widetilde{b}_{-4,1},\\[2pt]
[\![C_{0,1},\,C_{-4,0}]\!]&=-12\Omega_{-4,1}+\partial
b_{-4,1},\\[2pt]
[\![C_{-4,1},\,{C}_{-4,-4}]\!]&=30\,\Omega_{-4,1}-\widetilde{\Omega}_{-4,1}+\partial
\overline{b}_{-4,1}
\end{array}\end{align*}
where
\begin{align*}\begin{array}{llll}\widetilde{b}_{-4,1}(X)(f)&=-\frac{17}{13} X^{(5)}f'-\frac{30}{13}(X^{(4)}f''+X''f^{(4)})
-\frac{30}{13}X^{(3)}f^{(3)}+\frac{1}{5}X'f^{(5)},\\[2pt]
{b}_{-4,1}(X)(f)&=-\frac{17}{13}
X^{(5)}f'-\frac{30}{13}(X^{(4)}f''+X''f^{(4)})-\frac{30}{13}X^{(3)}f^{(3)},\\[2pt]
\overline{b}_{-4,1}(X)(f)&=\frac{36}{13}
X^{(5)}f'+\frac{495}{52}(X^{(4)}f''+X''f^{(4)})+\frac{345}{26}X^{(3)}f^{(3)}.
\end{array}\end{align*}
Thus, we obtain the integrability conditions corresponding to the
case $k=5$.

3) If $\lambda\notin\{a_1,a_2\}$   then $B_{\lambda,\lambda+6}$ is
coboundary. More precisily, $B_{\lambda,\lambda+6}$ is defined by
\begin{align*}
B_{\lambda,\lambda+6}=&t_{\lambda,\lambda+2}
t_{\lambda+2,\lambda+6}[\![C_{\lambda+2,\lambda+6},C_{\lambda,\lambda+2}]\!]+t_{\lambda,\lambda+3}
t_{\lambda+3,\lambda+6}[\![C_{\lambda+3,\lambda+6},C_{\lambda,\lambda+3}]\!]\\[2pt]&+t_{\lambda,\lambda+4}
t_{\lambda+4,\lambda+6}[\![C_{\lambda+4,\lambda+6},C_{\lambda,\lambda+4}]\!],
\end{align*}but, we show that
\begin{align*}\begin{array}{llll}
[\![C_{\lambda+2,\lambda+6},C_{\lambda,\lambda+2}]\!]=\partial
b_{\lambda,\lambda+6},\,\,
[\![C_{\lambda+3,\lambda+6},C_{\lambda,\lambda+3}]\!]=\partial
\widetilde{b}_{\lambda,\lambda+6}\,\,\text{and}\,\,[\![C_{\lambda+4,\lambda+6},C_{\lambda,\lambda+4}]\!]=\partial
\overline{b}_{\lambda,\lambda+6}
\end{array}\end{align*}
where
\begin{align*}\begin{array}{lllllllll}
b_{\lambda,\lambda+6}(X)(f)=&\frac{1}{14}\lambda(2\lambda^2+10\lambda+3)
\Big((-12-9\lambda+97\lambda^2+90\lambda^3+24\lambda^4)
 X^{(6)}f\\[2pt]
~&+(-72-404\lambda-41\lambda^2+127\lambda^3+60\lambda^4)X^{(5)}f''\\[2pt]
~&+(-180-163\lambda+83\lambda^2+160\lambda^3+80\lambda^4)X^{(4)}f^{(3)}\\[2pt]
~&+
(-240-922\lambda-13\lambda^2+155\lambda^3+60\lambda^4)X^{(3)}f^{(4)}\\[2pt]
~&+(-180-569\lambda-15\lambda^2+90\lambda^3+24\lambda^4)X''f^{(5)}\Big),\\[2pt]
\widetilde{b}_{\lambda,\lambda+6}(X)(f)=
&\frac{1}{14}\lambda(2\lambda^2+10\lambda+3)
\Big((-12\lambda^2-\lambda)
 X^{(6)}f-2(44\lambda^2+10\lambda)
X^{(5)}f''\\[2pt]
&~+(-68\lambda^2-32\lambda)X^{(4)}f^{(3)} +(12\lambda^2+4\lambda)
X''f^{(5)}\Big),\\[2pt]
\overline{b}_{\lambda,\lambda+6}(X)(f)=
&\frac{1}{14}\lambda(2\lambda^2+10\lambda+3)
\Big((-97\lambda^2-118\lambda^3) X^{(6)}f\\[2pt]
&~-(16\lambda+183\lambda^2+118\lambda^3)X^{(5)}f''-
(5\lambda+195\lambda^2+580\lambda^3)X^{(4)}f^{(3)}\\[2pt]
&~-(58\lambda+323\lambda^2+435\lambda^3)X^{(3)}f^{(4)}
-(103\lambda+377\lambda^2+174\lambda^3)X''f^{(5)}\Big).
\end{array}\end{align*}

For $\lambda=a_i$, the cohomology space $\mathrm{H}^2_{\rm
diff}(\mathfrak{vect}(1),{\mathcal{D}}_{a_i,a_i+6})$
is spanned by the cohomology classes of the 2-cocycles
$\Omega_{a_i,a_i+6}=[\![C_{a_i+3,a_i+6},\,{C}_{a_i,a_i+3}]\!]$ and
$\widetilde{\Omega}_{a_i,a_i+6}=[\![C_{a_i+6,a_i+6},\,{C}_{a_i,a_i+6}]\!]$.
Moreover, in this case, we have
\begin{align*}
B_{a_i,a_i+6}=&\,t_{a_i,a_i}
t_{a_i,a_i+6}[\![C_{a_i,a_i+6},C_{a_i,a_i}]\!]+{t}_{a_i,a_i+2}
t_{a_i+2,a_i+6}[\![C_{a_i+6,a_i+2},{C}_{a_i,a_i+2}]\!]\\[2pt]&+t_{a_i,a_i+3}
t_{a_i+3,a_i+6}[\![C_{a_i+3,a_i+6},C_{a_i,a_i+3}]\!]+t_{a_i,a_i+4}
t_{a_i+4,a_i+6}[\![C_{a_i+4,a_i+6},C_{a_i,a_i+4}]\!].
\end{align*}
But, by direct computation, we show that
\begin{align*}\begin{array}{llllllllll}
[\![C_{a_i,a_i+6},\,{C}_{a_i,a_i}]\!]&=-R_i\Omega_{a_i,a_i+6}-\widetilde{\Omega}_{a_i,a_i+6}+\partial
b_{a_i,a_i+6},\\[2pt]
[\![C_{a_i+2,a_i+6},\,{C}_{a_i,a_i+2}]\!]&=S_i\Omega_{a_i,a_i+6}+\partial
\widetilde{b}_{a_i,a_i+6},\\[2pt]
[\![C_{a_i+4,a_i+6},\,{C}_{a_i,a_i+4}]\!]&=-T_i\Omega_{a_i,a_i+6}+\partial
\overline{b}_{a_i,a_i+6}
\end{array}\end{align*} where
the maps $b_{a_i,a_i+6}$, $\widetilde{b}_{a_i,a_i+6}$ and
$\overline{b}_{a_i,a_i+6}$ are all proportional to the map $b$
defined by $$b(X)(f)=X^{(5)}f''.$$ We omit here the explicit
expressions of the scalar factors because they are  too
complicated. Thus, we obtain the integrability conditions
corresponding to the case $k=6$.
\end{proofname}
\begin{prop}\label{pr4}For $k=7,\,8,\,9,\,10$, we have the following
second-order integrability  conditions of the infinitesimal
deformation~(\ref{InDef}), where in the first line
$\lambda\notin\{0,-2, -4, -6\}$ and in the fourth
$\lambda\notin\{-7, 0,a_i, a_i-2,-4, -\frac{7\pm\sqrt{39}}{2}\}$:
\begin{equation}\label{k789}\begin{array}{llll}
(2\lambda+13)t_{\lambda,\lambda+3}t_{\lambda+3,\lambda+7}
+(1-2\lambda)t_{\lambda,\lambda+4}t_{\lambda+4,\lambda+7}&=&0,\\[2pt]
45~t_{-2,0}t_{0,5}-36~t_{-2,1}t_{1,5}-20~t_{-2,2}t_{2,5}&=&0,\\[2pt]
20~t_{-4,-1}t_{-1,3}+36~t_{-4,0}t_{0,3}+45~t_{-4,1}t_{1,3}&=&0,\\[2pt]
t_{\lambda,\lambda+4}t_{\lambda+4,\lambda+8}&=&0,\\[2pt]
60~t_{-7,-3}t_{-3,1}+~t_{-7,-4}t_{-4,1}&=&0,\\[2pt]
-60~t_{0,4}t_{4,8}+~t_{0,5}t_{5,8}&=&0,\\[2pt]
4~t_{-4,0}t_{0,4}-~t_{-4,1}t_{1,4}&=&0,\\[2pt]
\eta_i~t_{a_i-2,a_i}t_{a_i,a_i+6}+\theta_i~t_{a_i-2,a_i+2}t_{a_i+2,a_i+6}&=&0,\\[2pt]
\mu_it_{a_i,a_i+6}t_{a_i+6,a_i+8}+\nu_it_{a_i,a_i+4}t_{a+4,a_i+8}&=&0, \\[2pt]
t_{a_i,a_i+6}t_{a_i+6,a_i+9}&=&0,\\[2pt]
t_{a_i-3,a_i}t_{a_i,a_i+6}&=&0,\\[2pt]
t_{-8,-4}t_{-4,1}&=&0,\\[2pt]
t_{-4,0}t_{0,5}-t_{-4,1}t_{1,5}&=&0,\\[2pt]
t_{0,5}t_{5,9}&=&0,\\[2pt]
t_{a_i,a_i+6}t_{a_i+6,a_i+10}&=&0,  \\[2pt]
t_{a_i-4,a_i+}t_{a_i,a_i+6}&=&0
\end{array}\end{equation}
where
\begin{equation}\label{const3}\begin{array}{ll}\eta_1=7(76437+53739\sqrt{19}),
~~~\theta_1=64(1160123+30689\sqrt{19}),\\
\mu_1=8947638+205273\sqrt{19},~~
~\nu_1=96(474174+108783\sqrt{19})\end{array}\end{equation}
($\eta_2,~\theta_2,~\mu_2,~\nu_2$ are the conjugates respectively
of $\eta_1,~\theta_1,~\mu_1,~\nu_1$).

\end{prop}
\begin{proofname}. 1)  If $\lambda\notin\{ 0,-2,-4,-6\}$, we have
\begin{align*}
B_{\lambda,\lambda+7}=t_{\lambda,\lambda+3}
t_{\lambda+3,\lambda+7}[\![C_{\lambda+3,\lambda+7},C_{\lambda,\lambda+3}]\!]+t_{\lambda,\lambda+4}
t_{\lambda+4,\lambda+7}[\![C_{\lambda+4,\lambda+7},C_{\lambda,\lambda+4}]\!],
\end{align*} but we show that
\begin{align*}\begin{array}{llll}
[\![C_{\lambda+4,\lambda+7},C_{\lambda,\lambda+4}]\!]={1-2\lambda\over2\lambda+13}\Omega_{\lambda,\lambda+7}+\partial
b_{\lambda,\lambda+7}.
\end{array}\end{align*}
So, we obtain the second-order integrability conditions for $k=7$
and for generic $\lambda$. Besides, we study, as before, singular
values of $\lambda$ and then we obtain the corresponding
second-order integrability conditions. More precisely, the map
$B_{\lambda,\lambda+7}$ has the following form:
\begin{align*}
B_{\lambda,\lambda+7}~=\omega^1_{\lambda,\lambda+7}(t)
\Omega_{\lambda,\lambda+7}+\omega^2_{\lambda,\lambda+7}(t)\partial
b_{\lambda,\lambda+7},
\end{align*} where
\begin{align*}\omega^1_{\lambda,\lambda+7}(t)=
\left\{\begin{array}{llll}(2\lambda+13)t_{\lambda,\lambda+3}t_{\lambda+3,\lambda+7}
+(1-2\lambda)t_{\lambda,\lambda+4}t_{\lambda+4,\lambda+7}&
\hbox{ if }\,\,\lambda\notin\{ 0,-2,-4,-6\}\\[2pt]
-\frac{9}{4}~t_{-2,0}t_{0,5}+
\frac{9}{5}~t_{-2,1}t_{1,5}+~t_{-2,2}t_{2,5}
&\hbox{ if }\,\, \lambda= -2,\\[2pt]
\frac{5}{9}~t_{-4,-1}t_{-1,3}+~t_{-4,0}t_{0,3}+
\frac{5}{4}~t_{-4,1}t_{1,3}&\hbox{ if } \,\,\lambda=-4\\[2pt]
0&\hbox{ if } \,\,\lambda=0,-6.
\end{array}\right.\end{align*}

Hereafter we omit the expressions of the maps
$b_{\lambda,\lambda+k}$ and $\omega^2_{\lambda,\lambda+k}$ as they
are too long.

2) Now, for $k=8$ and $2\lambda\neq-7\pm\sqrt{39}$, the spaces
$\mathrm{H}^2_{\rm
diff}(\mathfrak{vect}(1),{\mathcal{D}}_{\lambda,\lambda+8})$
are spanned by the cohomology classes of the 2-cocycle
$\Omega_{\lambda,\lambda+8}=[\![C_{\lambda+4,\lambda+8},C_{\lambda,\lambda+4}]\!]$
and generically we have
\begin{align*}
B_{\lambda,\lambda+8}=t_{\lambda,\lambda+4}
t_{\lambda+4,\lambda+8}[\![C_{\lambda+4,\lambda+8},C_{\lambda,\lambda+4}]\!].
\end{align*}
But, for singular values of $\lambda$, other cup-products appear
in  the expression of $B_{\lambda,\lambda+8}$. More precisely, we
show that
\begin{align*}
B_{\lambda,\lambda+8}~=\omega^1_{\lambda,\lambda+8}(t)
\Omega_{\lambda,\lambda+8}+\omega^2_{\lambda,\lambda+8}(t)\partial
b_{\lambda,\lambda+8},
\end{align*} where
\begin{align*}\omega^1_{\lambda,\lambda+8}(t)=
\left\{\begin{array}{llll}t_{\lambda+4,\lambda+8} t_{\lambda,\lambda+4}&\hbox{
if }\,\,\lambda\neq 0,-7,-4,-\frac{7\pm\sqrt{39}}{2},a_i,a_i-2\\[2pt]
t_{-3,1}t_{-7,-3}+~\frac{1}{60}t_{-4,1}t_{-7,-4}&\hbox{ if }\,\, \lambda= -7,\\[2pt]
~t_{4,8} t_{0,4}-\frac{1}{60}~t_{0,5}t_{5,8}&\hbox{ if } \,\,\lambda=0\\[2pt]
t_{-4,0}t_{0,4}-\frac{1}{4}t_{-4,1}t_{1,4}&\hbox{ if }\,\, \lambda=-4\\[2pt]
\frac{\mu_i}{\nu_i}t_{a_i+6,a_i+8}t_{a_i,a_i+6}
+t_{a_i+4,a_i+8}t_{a_i,a_i+4}&\hbox{ if }\,\, \lambda=a_i\\[2pt]
\frac{\eta_i}{\theta_i}~t_{a_i,a_i+6}t_{a_i-2,a_i}
+~t_{a_i+2,a_i+6}t_{a_i-2,a_i+2}&\hbox{ if }\,\, \lambda=a_i-2\\[2pt]
0&\hbox{ if }\,\, \lambda=-\frac{7\pm\sqrt{39}}{2}.
\end{array}\right.\end{align*}

2) For $k=9$, the maps $B_{\lambda,\lambda+9}$ exist only for some
singular values of $\lambda$. More precisely, we show that, for
$\lambda\neq-4$,
\begin{align*}
B_{\lambda,\lambda+9}~=\omega_{\lambda,\lambda+9}(t)
\Omega_{\lambda,\lambda+9},
\end{align*}
where
\begin{align*}\omega_{\lambda,\lambda+9}(t)=\left\{\begin{array}{llll}0&\hbox{
if }\,\,\lambda\neq a_i,a_i-3,0,-8,-4,\\[2pt]
t_{-8,-4}t_{-4,1} &\hbox{ if } \,\,\lambda= -8,\\[2pt]
t_{0,5} t_{5,9}&\hbox{ if } \,\,\lambda=0,\\[2pt]
t_{a_i-3,a_i}t_{a_i,a_i+6} &\hbox{ if }\,\, \lambda=a_i-3,\\[2pt]
t_{a_i,a_i+6} t_{a_i+6,a_i+9}&\hbox{ if }\,\,
\lambda=a_i\end{array}\right.\end{align*} and
$$B_{-4,5}=(t_{-4,0}t_{0,5} -t_{-4,1}t_{1,5})\Omega_{-4,5}+t_{-4,1}t_{1,5}b_{-4,5}.$$

3) Finally, we show that
\begin{align*}
B_{\lambda,\lambda+10}~=\omega_{\lambda,\lambda+10}(t)
\Omega_{\lambda,\lambda+10},
\end{align*}
where
\begin{align*}\omega_{\lambda,\lambda+10}(t)=\left\{\begin{array}{llll}0&\hbox{
if }\,\,\lambda\neq a_i,a_i-4,\\[2pt]
t_{a_i,a_i+6}t_{a_i+6,a_i+10} &\hbox{ if }\,\, \lambda= a_i,\\[2pt]
t_{a_i-4,a_i}t_{a_i,a_i+6} &\hbox{ if }\,\, \lambda=a_i-4.\\[2pt]
\end{array}\right.\end{align*}\hfill$\Box$
\end{proofname}
 \vskip0.3cm

Now, in the following theorem, we recapitulate  the second-order
integrability conditions for the infinitesimal
deformation~(\ref{InfinDef2}). More precisely, we give the
necessary conditions to have the second term ${\cal L}^{(2)}$ of
(\ref{InfinDef2}). We give all conditions of second order, but,
any space ${\cal S}^n_{\delta}$ is concerned only by relations
between monomials
$t_{\lambda,\lambda+j}t_{\lambda+j,\lambda+k}{}$, where
$\delta-n\leq\lambda,\lambda+k\leq\delta$ and $0\leq j\leq k\leq10
$.

Our main result in this paper is the following

\begin{thm}
\label{Main12}  The conditions (\ref{k123}), (\ref{k456}) and
(\ref{k789}) are necessary and sufficient for second-order
integrability of the infinitesimal deformation~(\ref{InfinDef2}).
\end{thm}
\begin{proofname}. Of course, these conditions are necessary as it was shown
in Propositions \ref{pr1},  \ref{pr3} and  \ref{pr4}. Now, under
these conditions, the second term $\mathcal{L}^{(2)}$ of the the
$\frak{sl}(2)$-trivial infinitesimal deformation~(\ref{InfinDef2})
is a solution of the Maurer-Cartan equation (\ref{cap}). This
solution is defined up to a 1-coboundary and it has been shown
in~\cite{ff2,aalo2} that different choices of solutions of the
Maurer-Cartan equation correspond to equivalent deformations.
Thus, we can always choose
\begin{equation}\label{L2}
\begin{array}{llll}{\cal
L}^{(2)}=&\displaystyle\frac{1}{2}\sum_{\lambda}
t_{\lambda,\lambda+2} t_{\lambda+2,\lambda+4}b_{\lambda,\lambda+4}
+ \frac{1}{2}\widetilde{t}_{0,1} t_{1,4}\widetilde{b}_{0,4}+
\frac{1}{2}t_{-3,0}\widetilde{t}_{0,1} \widetilde{b}_{-3,1}
\\&+\displaystyle\frac{1}{2}\sum_{\lambda}t_{\lambda,\lambda+2}t_{\lambda+2,\lambda+5}b_{\lambda,\lambda+5}
+\displaystyle\frac{1}{2}\sum_{\lambda}t_{\lambda,\lambda+3}t_{\lambda+3,\lambda+5}\widetilde{b}_{\lambda,\lambda+5}\\
&+\displaystyle\frac{1}{2}\sum_{\lambda=0,-4}t_{\lambda,\lambda}t_{\lambda,\lambda+5}\overline{b}_{\lambda,\lambda+5}
+\displaystyle\frac{1}{2}\sum_{\lambda}t_{\lambda,\lambda+2}t_{\lambda+2,\lambda+6}b_{\lambda,\lambda+6}\\&
+\displaystyle\frac{1}{2}\sum_{\lambda}t_{\lambda,\lambda+3}t_{\lambda+3,\lambda+6}
\widetilde{b}_{\lambda,\lambda+6}+\frac{1}{2}\sum_{\lambda}t_{\lambda,\lambda+4}t_{\lambda+4,\lambda+6}
\overline{b}_{\lambda,\lambda+6}\\&+
\displaystyle\frac{1}{2}\sum_{\lambda}\omega^2_{\lambda,\lambda+7}(t)b_{\lambda,\lambda+7}+
\frac{1}{2}\sum_{\lambda}\omega^2_{\lambda,\lambda+8}(t)b_{\lambda,\lambda+8}+t_{-4,1}t_{1,5}b_{-4,5}.
\end{array}
\end{equation}
Of course, any $t_{\lambda,\lambda+k}$ appear in the expressions
of ${\cal L}^{(1)}$ or ${\cal L}^{(2)}$ if and only if
$\delta-\lambda$ and $k$ are integers satisfying
$\delta-n\leq\lambda,\lambda+k\leq\delta$. Theorem \ref{Main12} is
proved. \hfill$\Box$
\end{proofname}
{\rem There are no second-order conditions for integrability in
the following cases: \vskip0.3cm i) $k=0$, for all $\lambda,$

\vskip0.3cm ii) $k=1$ and $\lambda\neq0,$

\vskip0.3cm iii) $k=5$ and $\lambda\neq0,-4,$

\vskip0.3cm iv) $k=6$ and $\lambda\neq a_i,$

\vskip0.3cm v) $k=7$ and $\lambda=0,-6,$

 \vskip0.3cm vi) $k=8$ and
$\lambda=-\frac{7\pm\sqrt{39}}{2},$

\vskip0.3cm vii) $k=9$ and $\lambda\neq a_i,0,-4,-8,a_i-3,$

\vskip0.3cm viii) $k=10$ and $\lambda\neq a_i,a_i-4$.}

 \vskip0.3cm

\section{Examples}
The second-order conditions given in Theorem (\ref{Main12}) are
not, in general, sufficient, but they are in some cases. In this
section we give examples of symbol spaces ${\cal S}^n_{\delta}$
for which the corresponding second-order integrability conditions
are also sufficient and then we describe completely the formal
deformations of these spaces. Finally, we consider an example for
which the second-order integrability conditions are not
sufficient, but we exhibit the higher-order integrability
conditions and then we describe also, in this case, the formal
deformations.

{\ex \label{Example1} Consider ${\cal S}^2_{3}$. The infinitesimal
deformation of the $\mathfrak{vect}(1)$-action on ${\cal
S}^2_{3}$  is of the form $L_X+{\cal L}^{(1)}_X$, where $L_X$ is
the Lie derivative of ${\cal S}^2_{3}$  along the vector field
$X\frac{d}{dx}$ defined by(\ref{Lieder2}), and

\begin{equation}
\label{InfinDef3} {\cal
L}_X^{(1)}=\sum_{\lambda=1}^3\sum_{j=0}^{3-\lambda}{}t_{\lambda,\lambda+j}\,C_{\lambda,\lambda+j}(X)=
(t_{1,1}C_{1,1}+t_{1,3}C_{1,3}+t_{2,2}C_{2,2}+t_{3,3}C_{3,3})(X).
\end{equation}

We have the unique equation :
\begin{equation}\label{cond3}t_{1,3} (t_{1,1}-t_{3,3})=0\end{equation} as
necessary integrability condition of this infinitesimal
deformation. The following proposition shoes that this condition
is also sufficient.}

\begin{proposition} There are two deformations  of the
$\mathfrak{vect}(1)$-action on ${\cal S}^2_{3}$ with
three independent  parameters given by :

\begin{equation}
\label{ExDef1}{\cal L}_X=L_X+t_1{}(C_{1,1}(X)+
C_{3,3}(X))+t_2{}C_{2,2}(X)+t_3{}C_{1,3}(X),\end{equation} or

\begin{equation}
\label{ExDef2}{\cal L}_X=L_X+t_1{}C_{1,1}(X)
+t_2{}C_{2,2}(X)+t_3{}C_{3,3}(X).\end{equation}
\end{proposition}

\begin{proofname}. We consider the infinitesimal deformation
(\ref{InfinDef3}) and then we consider  solutions of
(\ref{cond3}). The first solution is : $t_{1,1}=t_{3,3}$, we put
$t_1=t_{1,1}=t_{3,3}$, $t_2=t_{2,2}$ and $t_3=t_{1,3}$ and then we
obtain (\ref{ExDef1}). The second solution is : $t_{1,3}=0$, we
put $t_1=t_{1,1}$, $t_2=t_{2,2}$ and $t_3=t_{3,3}$ and then we
obtain the second deformation. The solution $\mathcal{L}^{(2)}$ of
(\ref{MC1}) can be chosen identically zero. Choosing the
highest-order terms $\mathcal{L}^{(m)}$ with $m\geq3$, also
identically zero, one obviously obtains a deformation (which is of
order 1 in $t$). \hfill$\Box$

\end{proofname}

{\ex \label{Example2} Now consider ${\cal S}^3_{4}$. In this case
we have

\begin{equation}
\label{InfinDef4} \begin{array}{ll}{\cal L}_X^{(1)}
=&t_{1,1}C_{1,1}(X)+t_{1,3}C_{1,3}(X)+{t}_{1,4}{C}_{1,4}(X)+t_{2,2}C_{2,2}(X)\\[2pt]
&+t_{2,4}C_{2,4}(X)+t_{3,3}C_{3,3}(X)+t_{4,4}C_{4,4}(X).\end{array}
\end{equation}}
\begin{proposition}
The space ${\cal S}^3_{4}$ admits eight different formal
deformations with four independent parameters. They are all
equivalent to infinitesimal ones.
\end{proposition}

\begin{proofname}. The following equations are the necessary integrability
conditions of the infinitesimal deformation $L_X+{\cal L}_X^{(1)}$
(they are also sufficient):
\begin{equation}\label{cond2}\begin{array}{lll}t_{1,3} (t_{1,1}-t_{3,3})&=&0,\\[2pt]
t_{2,4} (t_{2,2}-t_{4,4})&=&0,\\[2pt]
t_{1,4}(t_{1,1}-t_{4,4})&=&0.\end{array}\end{equation}

There are eight solutions for these equations, so  ${\cal
S}^3_{4}$ admits eight different deformations with four
independent parameters. Like in the first example all these
deformations are equivalent to infinitesimal ones. \hfill$\Box$
\end{proofname}

{\ex \label{Example3} Consider ${\cal S}^3_{3}$.

\begin{equation}
\label{InfinDef5} \begin{array}{ll}{\cal L}_X^{(1)}
&=t_{0,0}C_{0,0}(X)+t_{0,1}C_{0,1}(X)+\widetilde{t}_{0,1}\widetilde{C}_{0,1}(X)+t_{0,2}C_{0,2}(X)
\\[2pt]
&+t_{0,3}C_{0,3}(X)+t_{1,1}C_{1,1}(X)+t_{2,2}C_{2,2}(X)+t_{3,3}C_{3,3}(X).\end{array}
\end{equation}

The integrability conditions are

\begin{equation}\label{cond4}\begin{array}{lll}
t_{0,1}(t_{0,0}-t_{1,1})&=&t_{1,1}{\widetilde t}_{0,1}\\[2pt]
t_{0,2} (t_{0,0}-t_{2,2})&=&0\\[2pt]
t_{0,3} (t_{0,0}-t_{3,3})&=&0.
\end{array}\end{equation}

In this case also these conditions are sufficient and any formal
deformation of ${\cal S}^3_{3}$ is equivalent to infinitesimal one
satisfying (\ref{cond4}).

One can construct a great number of examples of deformations of
${\cal S}^3_{3}$ with 4 (or less) independent parameters, but the
deformation $$ {\cal L}_X=L_X+{\cal L}^{(1)}_X
$$  is the miniversal one of ${\cal
S}^3_{3}$ with base $\mathcal{A}= \mathbb{C}[t]/\mathcal{R}$,
where $\mathbb{C}[t]=\mathbb{C}[t_{0,0},t_{0,1},\dots]$ and
$\mathcal{R}$ is the ideal generated by
$$t_{0,1}(t_{0,0}-t_{1,1})-t_{1,1}{\widetilde t}_{0,1},~~~~t_{0,2}
(t_{0,0}-t_{2,2})~\hbox{~ and ~}t_{0,3} (t_{0,0}-t_{3,3}).$$ }

{\ex \label{Example4} Consider ${\cal S}^4_{5}$. In this case we
have

\begin{equation}
\label{InfinDef6} \begin{array}{ll}{\cal L}_X^{(1)}
=&t_{1,1}C_{1,1}(X)+t_{1,3}C_{1,3}(X)+{t}_{1,4}{C}_{1,4}(X)+{t}_{1,5}{C}_{1,5}(X)\\[2pt]
&+t_{2,2}C_{2,2}(X)+t_{2,4}C_{2,4}(X)+{t}_{2,5}{C}_{2,5}(X)+t_{3,3}C_{3,3}(X)\\[2pt]
&+t_{3,5}C_{3,5}(X)+t_{4,4}C_{4,4}(X)+t_{5,5}C_{5,5}(X).\end{array}
\end{equation}}
\begin{proposition}
Any formal deformation of ${\cal S}^4_{5}$ is equivalent to a
polynomial one with degree $\leq3.$
\end{proposition}
\begin{proofname}. The second-order integrability conditions in this case are

\begin{equation}\label{cond5}\begin{array}{lll}t_{1,3} (t_{1,1}-t_{3,3})&=&0,\\[2pt]
t_{2,4} (t_{2,2}-t_{4,4})&=&0,\\[2pt]
t_{3,5} (t_{3,3}-t_{5,5})&=&0,\\[2pt]
t_{1,4}(t_{1,1}-t_{4,4})&=&0,\\[2pt]
t_{2,5}(t_{2,2}-t_{5,5})&=&0,\\[2pt]
t_{1,5} (t_{1,1}-t_{5,5})&=&0.\end{array}\end{equation} Under
these conditions the second-order term ${\cal
L}^{(2)}_X:\mathcal{F}_1\rightarrow\mathcal{F}_5$ is defined by
$$
\begin{array}{ll}{\cal
L}_X^{(2)}f=\frac{1}{2}\,t_{1,3}t_{3,5}b_{1,5}(X)f=-t_{1,3}t_{3,5}X^{(4)}f'.
\end{array}$$
The third-order term ${\cal
L}^{(3)}_X:\mathcal{F}_1\rightarrow\mathcal{F}_5$ is a solution of
the third-order Maurer-Cartan equation:
\begin{equation}\label{MC3}\partial{\cal L}^{(3)}=-\frac{1}{2}\,\sum_{i+j=3}[\![{\cal
L}^{(i)},{\cal L}^{(j)}]\!]. \end{equation} We compute the right
side of the equation (\ref{MC3}), so this equation becomes
$$
\begin{array}{ll}\partial{\cal L}^{(3)}(X,Y)f={1\over2}t_{1,1}t_{1,3}t_{3,5}X^{(4)}Y''f
+{1\over2}(t_{1,1}-t_{5,5})t_{1,3}
t_{3,5}X^{(4)}Y'f'-(X\leftrightarrow Y).
\end{array}
$$
Under the following third-order integrability condition:
\begin{equation}\label{cond6}
(t_{1,1}-t_{5,5})t_{1,3} t_{3,5}=0,\end{equation}
 the equation (\ref{MC3}) has a solution:
$$
\begin{array}{ll}{\cal L}_X^{(3)}f=\frac{1}{5}\,t_{1,1}t_{1,3}t_{3,5}X^{(5)}f.\end{array}
$$

Now, we compute the fourth-order term ${\cal
L}_X^{(4)}:\mathcal{F}_1\rightarrow\mathcal{F}_5$. It is a
solution of:

\begin{equation}\label{MC4}\partial{\cal L}^{(4)}=-\frac{1}{2}\,\sum_{i+j=4}[\![{\cal
L}^{(i)},{\cal L}^{(j)}]\!]=-\frac{1}{2}\,t_{1,1}[\![C_{1,1},{\cal
L}^{(3)}]\!]-\frac{1}{2}\,t_{5,5}[\![{\cal L}^{(3)},C_{5,5}]\!].
\end{equation}
It is easy to see that, under the conditions (\ref{cond5}) and
(\ref{cond6}), the right hand side of (\ref{MC4}) is identically
zero. Thus, the solution $\mathcal{L}^{(4)}$ of (\ref{MC4}) can be
chosen identically zero. Choosing the highest-order terms
$\mathcal{L}^{(m)}$ with $m\geq5$, also identically zero, one
obviously obtains a deformation (which is of order 3 in $t$).

\vskip.3cm Now, by studying  the equations (\ref{cond5}) and
(\ref{cond6}), we can see that, up to equivalence, the Lie
derivative on ${\cal S}^4_{5}$ admits a formal deformation with
seven independent parameters, this deformation corresponds to the
solution $t_{i,i}=t_{j,j}$ of the equations (\ref{cond5}) and
(\ref{cond6}). A great number of non trivial deformations with $k$
independent parameters can be constructed if $k<7$, each
deformation corresponds to a solution to equations (\ref{cond5})
and (\ref{cond6}). All these deformations are polynomial of order
equal or less than 3 in $t$.\hfill$\Box$\end{proofname}

 {\rem In the  previous four examples we obtain the same results
if we substitute ${\cal S}^m_{\lambda+n}$ for  ${\cal S}^m_{n}$
where $\lambda\in\mathbb{R}_+^*$. }
 \vskip0.3cm



\end{document}